\documentclass[reqno]{amsart}

\usepackage{caption} 
\usepackage{tikz}
\usepackage{wasysym} 
\usetikzlibrary{decorations.pathmorphing}
\usepackage{enumitem}
\def\bim{\begin{itemize}\item[]}
	\def\eim{\end{itemize}}

\usepackage{amsmath}
\usepackage{amssymb}
\allowdisplaybreaks
\usepackage{mathrsfs}

\usepackage{amsthm}
\newtheorem{theorem}{Theorem}

\def\[{[\! [}
\def\]{]\! ]}

\begin{document}

	%
	%
	
	\title[Horton-Strahler numbers for binary butterfly trees]
	{Horton-Strahler numbers for binary butterfly trees: exact analysis}

	\author[H. Prodinger ]{Helmut Prodinger }
	\address{Department of Mathematics, University of Stellenbosch 7602, Stellenbosch, South Africa
		and
		NITheCS (National Institute for
		Theoretical and Computational Sciences), South Africa.}
	\email{hproding@sun.ac.za, warrenham33@gmail.com}
	

	\begin{abstract}
		Peca suggested in a recent paper on the arxiv to consider binary butterfly trees and their Horton-Strahler numbers.
		The trees are obtained by glueing two binary trees together in a special way; the results are again binary trees but with 
		a different probability distribution. A thorough combinatorial analysis is provided and leads asymptotically to the same results as for classical binary trees.
		
	\end{abstract}
	
	
	\maketitle

	\section{Binary trees and Horton-Strahler numbers}
	\label{sec:Horton--Strahler}
	
	We start with some classical observations about binary trees and Horton-Strahler numbers (also called register function in Computer science). The pioneering papers are by Flajolet, Raoult, and Vuillemin \cite{FRV} and Kemp \cite{kemp}.
	The author has collected  some material   in~\cite{EATCS, collected-pap-intro, mansour}. 
	
	Binary trees may be expressed by the following symbolic equation, which says that they include the empty tree 
	and trees recursively built from a root followed by two subtrees (left and right), which are binary trees:
	\begin{center}
		\begin{tikzpicture}
			[inner sep=1.3mm,
			s1/.style={circle=10pt,draw=black!90,thick},
			s2/.style={rectangle,draw=black!50,thick},scale=0.5]
			
			\node at ( -4.8,0) { $\mathscr{B}$};
			
			\node at (-3,0) { $=$};
			\node(c) at (-1.5,0){ $\qed$};
			\node at (0.7,0) {$+$};
			\node(d) at (3,1)[s1]{};
			\node(e) at (2,-1){ $\mathscr{B}$};
			\node(f) at (4,-1){ $\mathscr{B}$};
			\path [draw,-,black!90] (d) -- (e) node{};
			\path [draw,-,black!90] (d) -- (f) node{};
			
		\end{tikzpicture}
	\end{center}
	
	Binary trees are counted by Catalan numbers and the  parameter \textsf{reg},  (register function, Horton-Strahler numbers) is recursively defined by attaching the number 0 to the leaves, and then working our way up: if both subtrees are labelled with the same number, the root will get 1 + this number, otherwise the larger value of the subtrees bubbles up. The value at the root is then the parameter of interest. 
	
	There is a recursive description of this function: $\textsf{reg}(\square)=0$, and if tree $t$ has subtrees $t_1$ and $t_2$, then
	\begin{equation*}
		\textsf{reg}(t)=
		\begin{cases}
			\max\{\textsf{reg}(t_1),\textsf{reg}(t_2)\}&\text{ if } \textsf{reg}(t_1)\ne\textsf{reg}(t_2),\\
			1+\textsf{reg}(t_1)&\text{ otherwise}.
		\end{cases}
	\end{equation*}
	
	Here is an example:
	
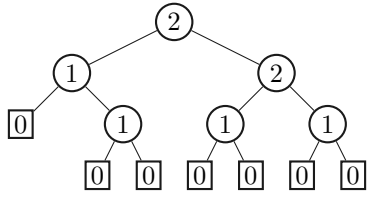
\begin{center}
		\begin{tikzpicture}
			[scale=0.34,inner sep=0.7mm,
			s1/.style={circle,draw=black!90,thick},
			s2/.style={rectangle,draw=black!90,thick}]
			\node(a) at ( 0,8) [s1] [text=black]{$2$};
			\node(b) at ( -4,6) [s1] [text=black]{$1$};
			\node(c) at ( 4,6) [s1] [text=black]{$2$};
			\node(d) at ( -6,4) [s2] [text=black]{$0$};
			\node(e) at ( -2,4) [s1] [text=black]{$1$};
			\node(f) at ( 2,4) [s1] [text=black]{$1$};
			\node(g) at ( 6,4) [s1] [text=black]{$1$};
			\node(h) at ( -3,2) [s2] [text=black]{$0$};
			\node(i) at ( -1,2) [s2] [text=black]{$0$};
			\node(j) at ( 1,2) [s2] [text=black]{$0$};
			\node(k) at ( 3,2) [s2] [text=black]{$0$};
			\node(l) at ( 5,2) [s2] [text=black]{$0$};
			\node(m) at ( 7,2) [s2] [text=black]{$0$};
			\path [draw,-,black!90] (a) -- (b) node{};
			\path [draw,-,black!90] (a) -- (c) node{};
			\path [draw,-,black!90] (b) -- (d) node{};
			\path [draw,-,black!90] (b) -- (e) node{};
			\path [draw,-,black!90] (c) -- (f) node{};
			\path [draw,-,black!90] (c) -- (g) node{};
			\path [draw,-,black!90] (e) -- (h) node{};
			\path [draw,-,black!90] (e) -- (i) node{};
			\path [draw,-,black!90] (f) -- (j) node{};
			\path [draw,-,black!90] (f) -- (k) node{};
			\path [draw,-,black!90] (g) -- (l) node{};
			\path [draw,-,black!90] (g) -- (m) node{};
		\end{tikzpicture}
	\end{center}
	
	Let $\mathscr{R}_{p}$ denote the family of binary
	trees with Horton-Strahler number equal to $p$, then one gets immediately from the recursive 
	definition:
 
	\begin{center}
		\begin{tikzpicture}
			[inner sep=1.3mm,
			s1/.style={circle=10pt,draw=black!90,thick},
			s2/.style={rectangle,draw=black!50,thick},scale=0.5]
			
			\node at ( -5,0) { $\mathscr{R}_p$};
			
			\node at (-4,0) { $=$};
			\node(a) at (-2,1)[s1]{};
			\node(b) at (-3,-1){ $\mathscr{R}_{p-1}$};
			\node(c) at (-1,-1){ $\mathscr{R}_{p-1}$};
			\path [draw,-,black!90] (a) -- (b) node{};
			\path [draw,-,black!90] (a) -- (c) node{};
			\node at (0.7,0) {$+$};
			\node(d) at (3,1)[s1]{};
			\node(e) at (2,-1){ $\mathscr{R}_{p}$};
			\node(f) at (4,-1.2){ $\sum\limits_{j<p}\mathscr{R}_{j} $};
			\path [draw,-,black!90] (d) -- (e) node{};
			\path [draw,-,black!90] (d) -- (f) node{};
			\node at (5+0.7,0) {$+$};
			\node(dd) at (5.5+3,1)[s1]{};
			\node(ee) at (5.5+2,-1.2){ $\sum\limits_{j<p}\mathscr{R}_{j}$};
			\node(ff) at (5.5+4,-1){ $\mathscr{R}_{p}$};
			\path [draw,-,black!90] (dd) -- (ee) node{};
			\path [draw,-,black!90] (dd) -- (ff) node{};
		\end{tikzpicture}
	\end{center}
	In terms of generating functions, these equations are translated into
	\begin{equation*}
		R_p(z)=zR_{p-1}^2(z)+2zR_p(z)\sum_{j<p}R_j(z);
	\end{equation*}
	the variable $z$ is used to mark the size (i.\,e., the number of internal nodes) of the binary tree.
	 See \cite{collected-pap-intro}; compare also~\cite{mansour}.
	 
	 Flajolet et al.\ resp.\ Kemp were able to solve this explicitly! Nowadays, several strategies to do this are known;
	 we only report the results as they will be used later in this paper.
	  The substitution
	\begin{equation*}
		z=\frac{u}{(1+u)^2}
	\end{equation*}
	that de Bruijn, Knuth, and Rice~\cite{BrKnRi72} also used, produces the nice expression
	\begin{equation*}
		R_p(z)=\frac{1-u^2}{u}\frac{u^{2^p}}{1-u^{2^{p+1}}}.
	\end{equation*}
	 The generating function $S_p(z)=R_{p}+R_{p+1}+R_{p+2}+\cdots$ of  binary trees with register function $\ge p$ is equally important. One can check directly that
	 \begin{center}
	 	\begin{tikzpicture}
	 		[inner sep=1.3mm,
	 		s1/.style={circle=10pt,draw=black!90,thick},
	 		s2/.style={rectangle,draw=black!50,thick},scale=0.5]
	 		
	 		\node at ( -5,0) { $\mathscr{S}_p$};
	 		
	 		\node at (-4,0) { $=$};
	 		\node(a) at (-2,1)[s1]{};
	 		\node(b) at (-3,-1){ $\mathscr{S}_{p-1}$};
	 		\node(c) at (-1,-1){ $\mathscr{S}_{p-1}$};
	 		\path [draw,-,black!90] (a) -- (b) node{};
	 		\path [draw,-,black!90] (a) -- (c) node{};
	 		\node at (0.7,0) {$+$};
	 		\node(d) at (3,1)[s1]{};
	 		\node(e) at (2,-1){ $\mathscr{S}_{p}$};
	 		\node(f) at (4.2,-1){ $\mathscr{B}\!\setminus \!S_{p-1} $};
	 		\path [draw,-,black!90] (d) -- (e) node{};
	 		\path [draw,-,black!90] (d) -- (f) node{};
	 		\node at (5+0.7,0) {$+$};
	 		\node(dd) at (5.5+3,1)[s1]{};
	 		\node(ee) at (5.3+2,-1.0){ $\mathscr{B}\!\setminus \!S_{p-1}$};
	 		\node(ff) at (5.7+4,-1){ $\mathscr{S}_{p}$};
	 		\path [draw,-,black!90] (dd) -- (ee) node{};
	 		\path [draw,-,black!90] (dd) -- (ff) node{};
	 	\end{tikzpicture}
	 \end{center}
	 Hence
	 \begin{align*}
	 	S_p(z)=zS_{p-1}^2+zS_p(B-S_{p-1})+ z(B-S_{p-1})S_p\quad p\ge1,\quad S_0 =B,
	 \end{align*}
	 and
	 \begin{equation*}
	 	S_p(z)=\frac{1-u^2}{u}\frac{u^{2^p}}{1-u^{2^{p}}}.
	 \end{equation*}
	 
	 Peca \cite{peca} introduced a ``butterfly'' tree by glueing two binary trees together. The motivation came from
	 binary \emph{search} trees which are binary trees but with a different statistics underlying; binary search trees are equivalent to permutations and thus enumerated by $n!$, not by Catalan numbers. The interest was, however, in the Horton-Strahler numbers of the
	 ``butterfly'' trees. The operation $\oplus$ is first defined on permutations: The operation $\pi_1\oplus \pi_2$ is defined like this:
	 Both permutations are written in one-line notation in terms of $1,2,\dots$; then, the second permutation is lifted up by $m$, if $m$ is the size of the permutation $\pi_1$; with the resulting string of numbers, the final binary search tree is formed. A thorough analysis of the Horton-Strahler numbers in these butterfly trees was requested in \cite{peca}. This is the purpose of this paper.
	 
	 The resulting tree can be described without permutations by finding the rightmost leaf of $t_1$, and then replacing it by the tree $t_2$. Peca also considered a version of $\ominus$, but that is basically a symmetric version of $\oplus$, and we will not  consider  it further. The resulting object is again a binary tree, but one cannot reconstruct the two binary trees $t_1$ and $t_2$ from it. In order to do so, we consider the path from the root to the rightmost leaf, and consider exactly one edge as being distinguished (marked). From this, the reconstruction is possible. We avoid to consider the empty tree  for $t_1$ since the path of interest has \emph{no} edges. 
	 
	 The generating function of these marked binary trees is $A:=(B-1)B=u(1+u)$. The enumeration is as follows, since $B$ is known:
	 \begin{equation*}
	 	A=\frac{1-3z}{2z^2}-\frac{(1-z)\sqrt{1-4z}}{2z^2},\qquad [z^n]A(z)=\frac{3(2n)!}{(n-1)!(n+2)!},\ n\ge1.
	 \end{equation*}

	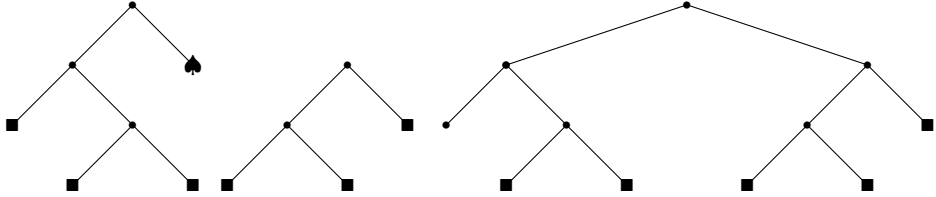
\begin{figure}
	 \begin{center}
	 \begin{tikzpicture}[scale=0.8]
	 	\path (-2,0) node(x1) {\tiny$\bullet$} ;
	 	\path (-1,1) node(x2) {\tiny$\bullet$};
	 	\path (-1,-1) node(x3) {\tiny$\bullet$};
	 	\path (-3,-1) node(x4) {\tiny$\blacksquare$}; 
	 	\path (0,0) node(x4) { \small$\spadesuit$}; 
	 	\path (-2,-2) node(x4) {\tiny$\blacksquare$}; 
	 	\path (0,-2) node(x4) {\tiny{$\blacksquare$}}; 
	 	\draw (-2,0) -- (-1,-1);
	 	\draw (-2,0) -- (-1,1) ;
	 	\draw (-3,-1) -- (-2,0) ;
	 	\draw (0,0) -- (-1,1);
	 	\draw (-2,-2) -- (-1,-1);
	 	\draw (0,-2) -- (-1,-1);
	 \end{tikzpicture}%
	 \begin{tikzpicture}[scale=0.8, xshift=4cm  ] 
	 \path (-2,0) node(x1) {\tiny$\bullet$} ;
	 \path (-1,1) node(x2) {\tiny$\bullet$};
	 \path (-1,-1) node(x3) {\tiny$\blacksquare$};
	 \path (-3,-1) node(x4) {\tiny$\blacksquare$}; 
	 \path (0,0) node(x4) {\tiny$\blacksquare$}; 
	 \draw (-2,0) -- (-1,-1);
	 \draw (-2,0) -- (-1,1) ;
	 \draw (-3,-1) -- (-2,0) ;
	 \draw (0,0) -- (-1,1);
	 \end{tikzpicture}
	 \begin{tikzpicture}[scale=0.8,xshift=8cm]
	 	\path (-2,0) node(x1) {\tiny$\bullet$} ;
	 	\path (1,1) node(x2) {\tiny$\bullet$};
	 	\path (-1,-1) node(x3) {\tiny$\bullet$};
	 	\path(-3,-1) node(x4) {\tiny$\bullet$}; 
	 	\path (4,0) node(x4) {\tiny$\bullet$}; 
	 	\path (-2,-2) node(x4) {\tiny$\blacksquare$}; 
	 	\path (0,-2) node(x4) {\tiny{$\blacksquare$}}; 
	 	\draw (-2,0) -- (-1,-1);
	 	\draw (-2,0) -- (1,1) ;
	 	\draw (-3,-1) -- (-2,0) ;
	 	\draw (4,0) -- (1,1);
	 	\draw (-2,-2) -- (-1,-1);
	 	\draw (0,-2) -- (-1,-1);
	 	\path (-2,0) node(x1) {\tiny$\bullet$} ;
 	 	\path (3,-1) node(x2) {\tiny$\bullet$};
 	 	\path (4,-2) node(x3) {\tiny$\blacksquare$};
 	 	\path ( 2,-2) node(x4) {\tiny$\blacksquare$}; 
 	 	\path (5,-1) node(x4) {\tiny$\blacksquare$}; 
 	 	\draw (4,0 ) -- (3 ,-1);
 	 	 \draw (4 ,-2) -- (3,-1) ;
 	 	\draw (2,-2) -- ( 3,-1);
 	 	 \draw ( 5,-1) -- (4,0);
	 \end{tikzpicture}
	 \end{center}
	 \caption{Two binary trees $t_1$, $t_2$, the rightmost leaf of $t_1$ indicated, and then $t_2$ glued there, with resulting binary tree
	 $t_1\oplus t_2$}
	\end{figure}
	
	\begin{figure}
		\begin{center}
			\begin{tikzpicture}[scale=0.8,xshift=8cm]
				\path (-2,0) node(x1) {\tiny$\bullet$} ;
				\path (1,1) node(x2) {\tiny$\bullet$};
				\path (-1,-1) node(x3) {\tiny$\bullet$};
				\path(-3,-1) node(x4) {\tiny$\bullet$}; 
				\path (4,0) node(x4) {\tiny$\bullet$}; 
				\path (-2,-2) node(x4) {\tiny$\blacksquare$}; 
				\path (0,-2) node(x4) {\tiny{$\blacksquare$}}; 
				\draw (-2,0) -- (-1,-1);
				\draw (-2,0) -- (1,1) ;
				\draw (-3,-1) -- (-2,0) ;
				\draw [red,thick](4,0) -- (1,1);
				\draw (-2,-2) -- (-1,-1);
				\draw (0,-2) -- (-1,-1);
				\path (-2,0) node(x1) {\tiny$\bullet$} ;
				\path (3,-1) node(x2) {\tiny$\bullet$};
				\path (4,-2) node(x3) {\tiny$\blacksquare$};
				\path ( 2,-2) node(x4) {\tiny$\blacksquare$}; 
				\path (5,-1) node(x4) {\tiny$\blacksquare$}; 
				\draw  (4,0 ) -- (3 ,-1);
				\draw (4 ,-2) -- (3,-1) ;
				\draw (2,-2) -- ( 3,-1);
				\draw ( 5,-1) -- (4,0);
			\end{tikzpicture}
			\begin{tikzpicture}[scale=0.6,xshift=8cm]
				\path (-2,0) node(x1) {\tiny$\bullet$} ;
				\path (1,1) node(x2) {\tiny$\bullet$};
				\path (-1,-1) node(x3) {\tiny$\bullet$};
				\path(-3,-1) node(x4) {\tiny$\bullet$}; 
				\path (4,0) node(x4) {\tiny$\bullet$}; 
				\path (-2,-2) node(x4) {\tiny$\blacksquare$}; 
				\path (0,-2) node(x4) {\tiny{$\blacksquare$}}; 
				\draw (-2,0) -- (-1,-1);
				\draw (-2,0) -- (1,1) ;
				\draw (-3,-1) -- (-2,0) ;
				\draw (4,0) -- (1,1);
				\draw (-2,-2) -- (-1,-1);
				\draw (0,-2) -- (-1,-1);
				\path (-2,0) node(x1) {\tiny$\bullet$} ;
				\path (3,-1) node(x2) {\tiny$\bullet$};
				\path (4,-2) node(x3) {\tiny$\blacksquare$};
				\path ( 2,-2) node(x4) {\tiny$\blacksquare$}; 
				\path (5,-1) node(x4) {\tiny$\blacksquare$}; 
				\draw (4,0 ) -- (3 ,-1);
				\draw (4 ,-2) -- (3,-1) ;
				\draw (2,-2) -- ( 3,-1);
				\draw [red,thick]( 5,-1) -- (4,0);
			\end{tikzpicture}
		\end{center}
		\caption{The resulting tree
			$t_1\oplus t_2$ can originate in two possible ways from glueing two trees together.}
	\end{figure}
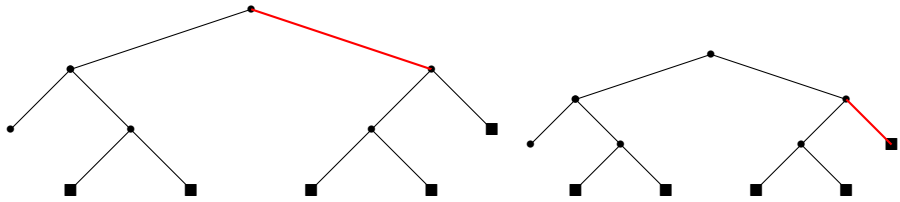
	When we consider the elements of $\mathscr{A}$ and ignore the marking of an edge on the rightmost path, we have a binary tree and thus can introduce the Horton-Strahler numbers on the elements of $\mathscr{A}$. No object has Horton-Strahler number equal to 0, and we introduce $T_p(z)$, the generating function of  the elements of $\mathscr{A}$ with Horton-Strahler number $\ge p$.
	We will find a recursion for $T_p$ that is akin  to the recursion for $S_p$ and use the fact that $S_p$ is explicitly known. This is similar to \cite{Yek} where the author solved a problem left open by Yekutieli and Mandelbrot~\cite{YM}, although the present situation is more involved. The further considerations deserve a section on its own.
	
	\section{Horton-Strahler numbers resulting from glueing together two binary trees}
	
	Recall that
	\begin{align*}
		S_p(z)=zS_{p-1}^2+zS_p(B-S_{p-1})+ z(B-S_{p-1})S_p,\quad p\ge1.
	\end{align*}
	Using a similar reasoning, distinguishing whether the first right edge is marked (first line) or is not marked (second line), we find
	\begin{align*}
		T_p(z)&=zS_{p-1}^2+zS_p(B-S_{p-1})+ z(B-S_{p-1})S_p\\
		&+zS_{p-1}T_{p-1}+zS_p(A-T_{p-1})+ z(B-S_{p-1})T_p,
	\end{align*}
	or
	 \begin{align*}
	 	T_p(z)&=S_p+zS_{p-1}T_{p-1}+zS_p(A-T_{p-1})+ z(B-S_{p-1})T_p,\quad p\ge2, \text{ and }T_1=A.
	 \end{align*}
	 Further
	 \begin{align*}
	 	T_p(z)\big(1-z(B-S_{p})\big)&=S_p(1+zA)+zT_{p-1}(S_{p-1}-S_p).
	 \end{align*}
	 A direct computation yields
	 \begin{align*}
	 	1-z(B-S_{p})
	 	&=\frac{1-u^{2^p+1}}{(1+u)(1-u^{2^p})} ,\quad 1+zA= \frac{1+u+u^2}{1+u}  .
	 \end{align*}
	 Further 
	 \begin{align*}
	 	 S_{p-1}-S_p=R_{p-1}=\frac{1-u^2}{u}\frac{u^{2^{p-1}}}{1-u^{2^{p}}}.
	 \end{align*}
	 Some rearrangements lead to
	 	 \begin{align*}
	 	 	T_p(z) (1-u^{2^{p-1}+1})(1+u^{2^{p-1}})   &=\frac{1-u^2}{u}u^{2^{p}}  (1+u+u^2) +T_{p-1}  (1-u)u^{2^{p-1}};
	 	 \end{align*}
	 	 	 this first order recursion will be solved by introduction of summation factors and eventually by summation. This will be done in a few steps;	 
	  \begin{align*}
	  	\frac{T_p(z)}{u^{2^p}} (1-u^{2^{p-1}+1})(1+u^{2^{p-1}})   &=\frac{1-u^2}{u}   (1+u+u^2) +\frac{T_{p-1}(z)}{u^{2^{p-1}}}  (1-u); 
	  \end{align*}
  \begin{align*}
	 	\frac{T_p(z)}{u^{2^p}} (1-u^{2^{p-1}+1})\prod_{j=0}^{p-1}(1+u^{2^{j}})   &=\frac{1-u^2}{u}   (1+u+u^2)\prod_{j=0}^{p-2}(1+u^{2^{j}})\\& +\frac{T_{p-1}(z)}{u^{2^{p-1}}}  (1-u) \prod_{j=0}^{p-2}(1+u^{2^{j}}).
	 \end{align*}
	 Evaluating the products yields
	 \begin{align*}
	 	\frac{T_p(z)(1-u^{2^{p-1}+1})  }{u^{2^p}} (1-u^{2^{p}})  &=\frac{1-u^2}{u}   (1+u+u^2)(1-u^{2^{p-1}}) \\&+\frac{T_{p-1}(z)(1-u^{2^{p-1}})}{u^{2^{p-1}}}  (1-u) 
	 \end{align*}
	 and	  
	\begin{align*}
		\frac{T_p(z)(1-u^{2^{p}})  }{u^{2^p}(1-u)^p} (1-u^{2^{p-1}+1})  &=\frac{1-u^2}{u(1-u)^p}   (1+u+u^2)(1-u^{2^{p-1}}) \\&+\frac{T_{p-1}(z)(1-u^{2^{p-1}})}{u^{2^{p-1}}(1-u)^{p-1}}   
	\end{align*} 
	Introducing a further product, 	
	\begin{align*}
		\frac{T_p(z)(1-u^{2^{p}})  }{u^{2^p}(1-u)^p}\prod_{j=0}^{p-1} (1-u^{2^{j}+1})  &=\frac{1+u}{u(1-u)^{p-1}}   (1+u+u^2)(1-u^{2^{p-1}})\prod_{j=0}^{p-2} (1-u^{2^{j}+1})\\& +\frac{T_{p-1}(z)(1-u^{2^{p-1}})}{u^{2^{p-1}}(1-u)^{p-1}}   \prod_{j=0}^{p-2} (1-u^{2^{j}+1}).
	\end{align*} 
	An abbreviation is useful:
	\begin{align*}
		\omega_p:=\frac{T_p(z)(1-u^{2^{p}})  }{u^{2^p}(1-u)^p} \prod_{j=0}^{p-1} (1-u^{2^{j}+1}),  
	\end{align*} 
	then we get a form that can be summed:
	\begin{align*}
		\omega_p=\omega_{p-1}+\frac{1+u}{u(1-u)^{p-1}}   (1+u+u^2)(1-u^{2^{p-1}})\prod_{j=0}^{p-2} (1-u^{2^{j}+1}) 
	\end{align*} 
	  and so	 
	\begin{align*}
		\omega_p=\omega_{1}+\sum_{h=1}^{p-1}\frac{1+u}{u(1-u)^{h}}   (1+u+u^2)(1-u^{2^{h}})\prod_{j=0}^{h-1} (1-u^{2^{j}+1}). 
	\end{align*} 
	 Coming back to the original quantities $T_p(z)$,	 
	\begin{align*}
		T_p(z)&=\omega_p\frac{u^{2^p}(1-u)^p }  { (1-u^{2^{p}})  }\bigg/\prod_{j=0}^{p-1} (1-u^{2^{j}+1})\\	 
		 &=\bigg[(1+u)^2   +(1+u+u^2)\sum_{h=1}^{p-1}   \frac{1-u^{2^{h}}}{1-u}\prod_{j=0}^{h-1} \frac{1-u^{2^{j}+1}}{1-u}\bigg] \\
		&\qquad\times\frac{1-u^2}{u}
		\frac{u^{2^p}  }  {  1-u^{2^{p}}   } \prod_{j=0}^{p-1}\frac{1-u} {1-u^{2^{j}+1}}.
	\end{align*}
	Note that $\dfrac{1-u^2}{u}\dfrac{u^{2^p}  }  {  1-u^{2^{p}}   } =S_p(z)$.

	\begin{theorem} The generating function $T_p(z)$ of trees in $\mathscr{A}$ with Horton-Strahler number $\ge p$ is for $p\ge0$
		 given by
		 \begin{align*}
		 	T_p(z) =S_p(z)\bigg[(1+u)^2   +(1+u+u^2)\sum_{h=1}^{p-1}   \frac{1-u^{2^{h}}}{1-u}\prod_{j=0}^{h-1} \frac{1-u^{2^{j}+1}}{1-u}\bigg] 
		 	 \prod_{j=0}^{p-1}\frac{1-u} {1-u^{2^{j}+1}}.
		 \end{align*}
		 
		\end{theorem}
		Note that this formula has been tested. Since the generating function is fully explicit one has (in principle) access to the coefficients, i.\,e., to the numbers of binary butterfly trees of a given number of nodes and a given Horton-Strahler number.
	
	\clearpage
	
		\section{The average value of Horton-Strahler numbers in marked binary trees}

	  By general principles, the generating function
	 \begin{align*}
	 \sum_{p\ge1}T_p(z)
	 \end{align*} 
	is, apart from normalization, the generating function of the averages; note that it was
	\begin{align*}
		\sum_{p\ge1}S_p(z)
	\end{align*}
	in the classical case, and the latter series is well-understood.
	
	To understand the strategy, we need to expand the generating function around $u=1$; in some instances it helps to
	set $u=e^{-t}$ and expand around $t=0$. Consider
	 	\begin{align}\label{ugly}
		\bigg[(1+u)^2   +(1+u+u^2)\sum_{h=1}^{p-1}   \frac{1-u^{2^h}}{1-u}\prod_{j=0}^{h-1} \frac{1-u^{2^{j}+1}}{1-u}\bigg] 
		\prod_{j=0}^{p-1}\frac{1-u} {1-u^{2^{j}+1}};
			\end{align}
			the term $\dfrac{1-u^2}{u}\displaystyle{\sum\limits_{p\ge1}\frac{u^{2^p}}{1-u^{2^{p}}}}$   will be brought in later.
	
	The ugly term \eqref{ugly} is actually simpler to handle, since we are only interested in the leading term, which is a constant (in the expansion around $t=0$). For the following we might replace a factor $\frac{1-u^{d}}{1-u}$ by $d$, and instead of \eqref{ugly}
			consider
			\begin{align}\label{ugly2}
				\lambda_p:=\bigg[4   +3\sum_{h=1}^{p-1}    2^h\prod_{j=0}^{h-1} (2^{j}+1)\bigg] 
				\prod_{j=0}^{p-1}\frac{1 } { 2^{j}+1}.
			\end{align}
			The sequence $\lambda_p$ converges to 3 exponentially fast.
			We consider now
			\begin{align*} 
				\sum_{p\ge1}\frac{u^{2^p}}{1-u^{2^p}}\lambda_p =\sum_{p,k\ge1}\lambda_p e^{-t{k2^p}}.
			\end{align*}
			The Mellin transformation will be applied, as in many related projects:
			\begin{align*} 
				\mathscr{M} \sum_{p,k\ge1}\lambda_p e^{-tk2^p}=\zeta(s)\mathscr{M} \sum_{p\ge1}\lambda_p e^{-t2^p}
				=\zeta(s)\Gamma(s) \underbrace{\sum_{p\ge1}\lambda_p  2^{-ps}}_{ \Lambda(s):=}.
			\end{align*}
			Now
			\begin{align*} 
				 \Lambda(s)=\sum_{p\ge1}\lambda_p  2^{-ps}=\frac{3}{2^s-1}+ \sum_{p\ge1}\mathscr{O}(2^{-p})  2^{-ps}.
			\end{align*}
			The remainder has shifted singularities, and thus the dominant ones are at $\Re s=0$. The general plan is to look at the
			 inverse Mellin transform $\zeta(s)\Gamma(s)  \Lambda(s)t^{-s}$   and   at its residues. Since the dominant ones are
			at $\Re s=0$, we can concentrate on $3\zeta(s)\Gamma(s)\frac{t^{-s}}{2^s-1}$, which is just 3 times the relevant quantity for classical binary trees.
			
			For the final averages of the Horton-Strahler numbers, we have to divide by
			\begin{align*} 
				[z^n]A(z)=\frac{3(2n)!}{(n-1)!(n+2)!}=\frac{3(2n)!}{n!(n+1)!}\Big (1+\mathscr{O}\Big(\frac1n\Big)\Big),
			\end{align*}
			so that within the accuracy the average Horton-Strahler numbers are computed, the asymptotic formula is the same:
			\begin{theorem}
				The average value of the Horton-Strahler numbers (register function) of all butterfly trees in the sense of Peca of size $n$
			 is given by the asymptotic formula (with $\chi_k=\frac{2k\pi i}{\log2}$)
				\begin{align*}
					\log_4 n&- \frac{\gamma }{2\log2}- \frac{1}{\log2}+\frac12+ \log_2 \pi   
					+\frac{1 }{\log2}\sum_{k\neq0} \zeta(\chi_k) \Gamma\big(\frac{\chi_k}{2}\big)(\chi_k-1)   n^{\chi_k/2}\\
					&=\log_4 n- \frac{\gamma }{2\log2}- \frac{1}{\log2}+\frac12+\log_2 \pi +\psi(\log_4n),
				\end{align*}
				with a tiny periodic function $\psi(x)$ of period 1.
			\end{theorem}	
			These oscillations are usually bounded by $10^{-5}$, say. See~\cite{FRV} for some explicit error bounds in the classical case.
			  The remainder term in the asymptotic formula is of the form $\mathscr{O}((\log^\ast\! n)/n)$ and has never been computed explicitly neither in the classical case nor for the butterfly trees because of the complexity of the computations.
			  
			  \bibliographystyle{plain} 

	\end{document}